# The Bi-Musquash Conjecture

by Jon Henry Sanders
jon_sanders@partech.com

Bi-Musquash Conjecture

A construction (attributed to Conway) for 'doubling' a thrackled m-gon of odd length m to get a thrackled 2m-gon (of even length) was given in [1]. Has anyone noticed that if this construction is applied to an odd length musquash we get a 'bi-musquash', i.e., a thrackled $n$-gon, $n = 2m$, m odd, with the property that row 1 of the intersection table generates the odd rows and row 2 generates the even rows? More precisely, for the definition of a 'bi-musquash' we require that[1] :

(d') If $e_1$ intersects edges in the order: $e_{k_1}, \ldots, e_{k_{n-3}}$, and $e_2$ intersects edges in the order $e_{j_1}, \ldots, e_{j_{n-3}}$ then for (odd) $i = 3, 5, \ldots, n-1$, edge $e_i$ intersects edges in the order $e_{k_1+i-1}, \ldots, e_{k_{n-3}+i-1}$, and for (even) $i = 4, \ldots, n$, edge $e_i$ intersects edges in the order : $e_{j_1+i-2}, \ldots, e_{j_{n-3}+i-2}$, where the edge subscripts are computed modulo $n$. Note that by this definition, an even musquash is automatically a bi-musquash. Remarks of [2,3] imply that bi-musquashes must have n/2 rotational symmetry and indeed the drawing of the bi-musquash of length 14 [Fig.1], as well the 'chordal representation'[2] of it's Gauss word [Fig.2] exhibit this. I guess the conjecture to prove[3] would be that bi- musquashes do not exist for n = 4m and that they are essentially[4] unique for each n, with n = 2 (mod 4). The intersection table of bi-musquash 14 is shown below.

| 8  | 7  | 12 | 5  | 10 | 3  | 13 | 6  | 11 | 4  | 9  |
|----|----|----|----|----|----|----|----|----|----|----|
| 9  | 8  | 6  | 13 | 4  | 11 | 7  | 14 | 5  | 12 | 10 |
| 10 | 9  | 14 | 7  | 12 | 5  | 1  | 8  | 13 | 6  | 11 |
| 11 | 10 | 8  | 1  | 6  | 13 | 9  | 2  | 7  | 14 | 12 |
| 12 | 11 | 2  | 9  | 14 | 7  | 3  | 10 | 1  | 8  | 13 |
| 13 | 12 | 10 | 3  | 8  | 1  | 11 | 4  | 9  | 2  | 14 |
| 14 | 13 | 4  | 11 | 2  | 9  | 5  | 12 | 3  | 10 | 1  |
| 1  | 14 | 12 | 5  | 10 | 3  | 13 | 6  | 11 | 4  | 2  |
| 2  | 1  | 6  | 13 | 4  | 11 | 7  | 14 | 5  | 12 | 3  |
| 3  | 2  | 14 | 7  | 12 | 5  | 1  | 8  | 13 | 6  | 4  |
| 4  | 3  | 8  | 1  | 6  | 13 | 9  | 2  | 7  | 14 | 5  |
| 5  | 4  | 2  | 9  | 14 | 7  | 3  | 10 | 1  | 8  | 6  |
| 6  | 5  | 10 | 3  | 8  | 1  | 11 | 4  | 9  | 2  | 7  |
| 7  | 6  | 4  | 11 | 2  | 9  | 5  | 12 | 3  | 10 | 8  |

Intersection Table of Bi-Musquash 14

The first two rows generate the table in the required way. In general, the conjecture is that the first two rows (generators) of a bi-musquash of length 2p, p odd, are either:

SET1
**Row 1:** p+1, p, 2p-2, p-2 , 2p-4, p-4 , ..., p+3, 3, 2p-1, p-1, 2p-3, p-3, ..., p+4, 4, p+2
**Row 2:** p+2, p+1, p-1, 2p-1, p-3, 2p-3, ..., 4, p+4, p, 2p, p-2, 2p-2, ..., 5, p+5, p+3

or

SET2
**Row 1:** p, 2p-2, p-2 , 2p-4, p-4 , ..., p+3, 3, 2p-1, p-1, 2p-3, p-3, ..., p+4, 4, p+2, p+1
**Row 2:** p+1, p-1, 2p-1, p-3, 2p-3, ..., 4, p+4, p, 2p, p-2, 2p-2, ..., 5, p+5, p+3, p+2,

each giving rise to a mirror image of the other.

---

[1] This is a modification of condition (d) of [2].
[2] i.e., the arrangement of the symbols of a Gauss word around a circle with a chord drawn between each pair of identical symbols
[3] G. Cairns has pointed out to me that " there are two inequivalent musquash of the 6-gon; one being the Conway doubling of the 3-gon" . Thus the case n = 6 of the bi-musquash conjecture needs to be modified appropriately.
[4] i.e., there exist precisely two and they are reflections of each other

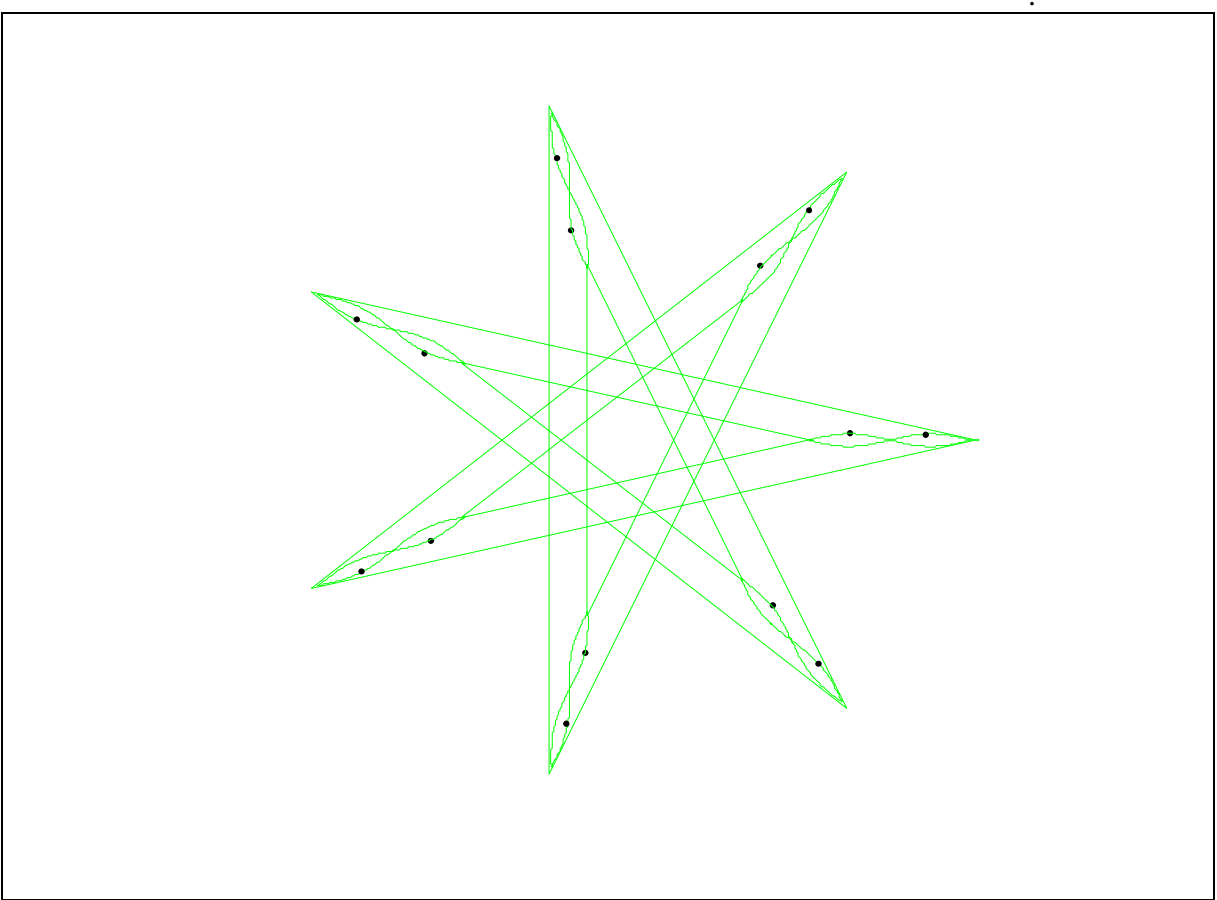

Figure 1  Bi-Musquash 14

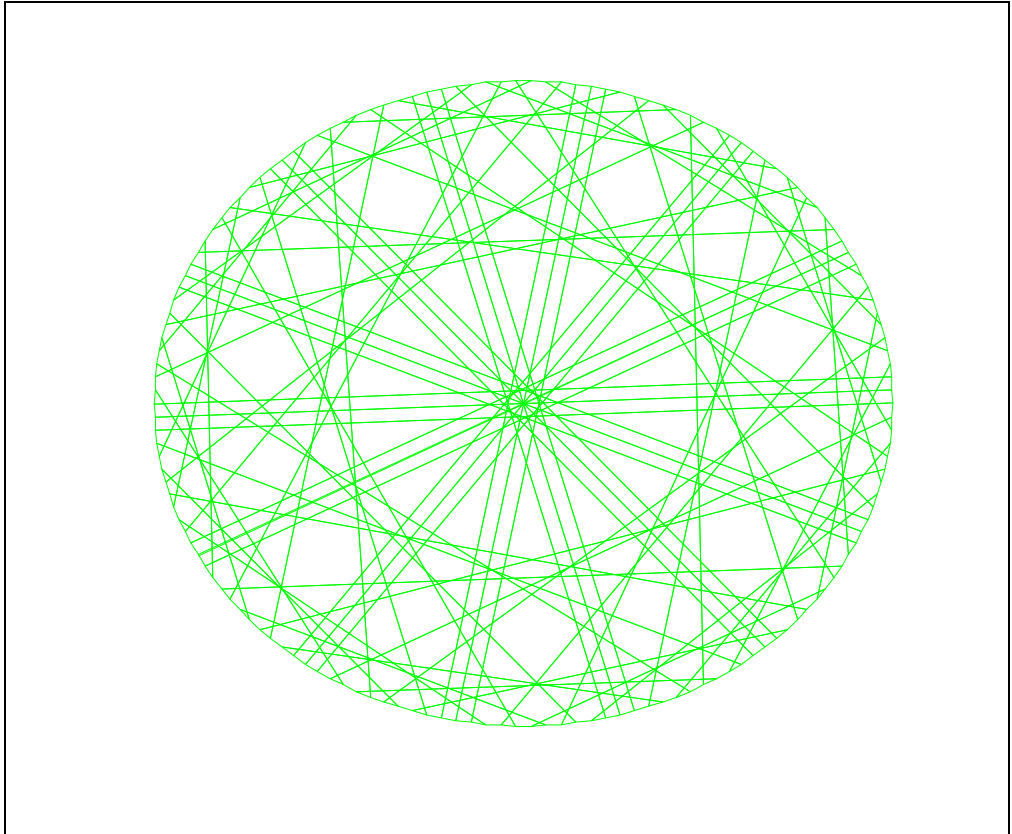

Figure 2  Chordal Representation of Bi-Musquash 14